\documentclass[12pt]{article}
\usepackage{amssymb,amsmath,amsthm}

\begin{document}

\title{An asymptotic property of quaternary additive codes}
\author{J\"urgen Bierbrauer\\
Department of Mathematical Sciences\\
Michigan Technological University\\
Houghton, Michigan 49931 (USA) \\ \\
S. Marcugini and F. Pambianco \\
Dipartimento di Matematica e Informatica\\
Universit\`a degli Studi di Perugia\\
Perugia (Italy) 
\footnote{The research of S. Marcugini and F. Pambianco was supported in part by the Italian National Group for Algebraic and Geometric Structures and their Applications (GNSAGA - INDAM) and by University of Perugia (Project No. 98751: Strutture Geometriche, Combinatoria e loro Applicazioni, Base Research Fund 2017-2019; Fighting Cybercrime with OSINT, Research Fund 2021). 
}}

\date{}

\maketitle
\newtheorem{Theorem}{Theorem}
\newtheorem{Proposition}{Proposition}
\newtheorem{Lemma}{Lemma}
\newtheorem{Definition}{Definition}
\newtheorem{Corollary}{Corollary}
\newtheorem{Example}{Example}
\def\nz{\mathbb{N}}
\def\gz{\mathbb{Z}}
\def\rz{\mathbb{R}}
\def\ef{\mathbb{F}}
\def\CC{\mathbb{C}}
\def\o{\omega}
\def\p{\overline{\omega}}
\def\e{\epsilon}
\def\a{\alpha}
\def\b{\beta}
\def\g{\gamma}
\def\d{\delta}
\def\l{\lambda}
\def\s{\sigma}
\def\bsl{\backslash}
\def\la{\longrightarrow}
\def\arr{\rightarrow}
\def\ov{\overline}
\def\sm{\setminus}
\newcommand{\D}{\displaystyle}
\newcommand{\T}{\textstyle}

\begin{abstract}
Let $n_k(s)$ be the maximal length $n$ such that a quaternary additive $[n,k,n-s]_4$-code exists.
We solve a natural asymptotic problem by determining the lim sup $\lambda_k$ of $n_k(s)/s,$ and the smallest value of $s$
such that $n_k(s)/s=\lambda_k.$ Our new family of quaternary additive codes has parameters
$[4^k-1,k,4^k-4^{k-1}]_4=[2^{2k}-1,k,3\cdot 2^{2k-2}]_4$ (where $k=l/2$ and $l$ is an odd integer).
These are constant-weight codes. The binary codes obtained by concatenation meet the Griesmer bound with equality.
The proof is in terms of multisets of lines in $PG(l-1,2).$
\end{abstract}

\section{Introduction}
\label{introsection}
In geometric terminology the theory of {\bf linear} codes is the theory of point sets in Galois spaces (projective geometries
over finite fields). For an introduction to geometries over finite fields and their connections with coding theory, see  
\cite{book2nded,ES2016,Giulietti2013,Hirschfeld1985,Hirschfeld1999,HirschfeldStorme,LandgevStorme}.
{\bf Additive} codes are a far-reaching and natural generalization of linear codes, see also Chapter 18 of \cite{book2nded}. 
 The quaternary case is of special interest, among others because of a close link to the
theory of quantum stabilizer codes and their geometric representations, see~\cite{quantgeom,CRRS,DHY}.
We study additive quaternary codes using geometric language. 
By definition such a code is a vector space over $\ef_2.$ If the vector space dimension is $l,$ then we call $k=l/2$ the
dimension of the code. It follows that dimensions $k$ of quaternary additive codes need not be integer
(but $l=2k$ is integer). Geometrically
the theory of additive quaternary codes is the study of multisets of lines in binary projective spaces.
An additive $[n,k,d]_4$-code is equivalent to a multiset of $n$ lines (the {\bf codelines}) in $PG(2k-1,2)$ such that each hyperplane contains at most $s=n-d$ of the codelines (in the multiset sense). For example, a code
$[7,1.5,6]_4$ corresponds to the set of lines of the Fano plane, and a code with parameters $[31,2.5,24]_4$ is equivalent to a
multiset of $31$ lines in $PG(4,2)$ such that each hyperplane $PG(3,2)$ contains at most $s=31-24=7$ of those
codelines. We prefer to work with the {\bf species} $s=n-d$ instead of the minimum distance.
 
\begin{Definition}
\label{optidef}
Given $s$ and a dimension $k\geq 1.5$ let $n_{k}(s)$ be the maximal length $n$ such that an additive quaternary
$[n,k,n-s]_4$-code exists.
\end{Definition}

The determination of the numbers $n_k(s)$ is a natural coding-theoretic optimization problem.
If $[n,k,d]_4$ exists then so does its {\bf concatenated code,} a binary linear $[3n,2k,2d]_2$-code.
The Griesmer bound for the concatenated code yields an existence condition for the additive code.
In earlier work we studied low dimensions. In \cite{dim2pt5} we determined the numbers $n_k(s)$ for all
dimensions $k\leq 3$ and all $s.$ In particular we showed that $[n,2.5,d]_4$ exists if and only if the parameters of the concatenated code $[3n,5,2d]_2$ satisfy the Griesmer bound.
Important recent contributions to the problem of determining the numbers $n_k(s)$ for small dimension $k$ 
are \cite{GLLL17,SLKS23}.\\
Our main result will be proved in Sections \ref{asympsection} and \ref{shortsection}.
In Section \ref{asympsection} we determine the maximal value $\lambda_k$ of $n/s$ for quaternary additive 
$[n,k,n-s]_4$ for each given dimension $k$ (Definition \ref{asympdef} and Theorem \ref{asympGriesmertheorem}).
In Section \ref{shortsection} we restrict to codes satisfying $n/s=\lambda_k$ and determine the minimal
length $n$ (equivalently: the minimal value of $s$) for which this is the case (Definition \ref{minsdef} and
Theorem \ref{asympGriesmertheorem}). Our main result is then the construction of a family of
$k$-dimensional quaternary additive codes with parameters $[4^k-1,k,4^k-4^{k-1}]_4$ for each
non-integer dimension $k.$ These are the $k$-dimensional quaternary additive codes of smallest
length $n$ among those with maximal $n/s.$ They are constant-weight codes.
Geometrically our family of codes corresponds to $PG(k-1,4)$ if $k$ is integer, to a $3$-cover of 
$PG(2k-1,2)$ if $k$ is non-integer (see Definition \ref{coverdef}, Lemma \ref{coverlemma}
and Theorem \ref{3covertheorem}).\\
In the final Section \ref{variantsection} we use a variant of the construction to obtain an interesting
infinite family of additive codes.

\section{An asymptotic problem}
\label{asympsection}

In the present work we study a natural and basic asymptotic version of the problem in arbitrary dimension:

\begin{Definition}
\label{asympdef}
For dimension $k\geq 1.5$ let $\lambda_k$ be the lim sup of $n_{k}(s)/s.$ 
\end{Definition}

\begin{Lemma}[sum construction]
\label{sumlemma}
If there is an $[n,k,d]_4$-code and an $[m,k,t]_4$-code, then there is an $[n+m,k,d+t]_4$-code.
\end{Lemma}
\begin{proof}
This corresponds to taking the union of the line sets, in the multiset sense.
\end{proof}

In particular $n_k(s_1+s_2)\geq n_k(s_1)+n_k(s_2).$ The sum construction also shows 
that for each $s$ the number $n_{k}(s)/s$ is a lower bound on $\lambda_{k}.$ \\
Let $g(l)=(2^l-1)(2^l-2)/6$ be the number of lines in $PG(l-1,2).$

\begin{Theorem}
\label{asympGriesmertheorem}
We have $\lambda_k=4+3/(4^{k-1}-1)=(4^k-1)/(4^{k-1}-1)$ for all $k\geq 1.5.$
For $s=g(2k-1)$ we have $n_k(s)/s=\lambda_k.$
\end{Theorem}
\begin{proof}
The Griesmer bound of the concatenated code yields an upper bound on $\lambda_k:$ 
if our additive code has parameters $[n,k,d]_4,$ then the concatenated code is a binary linear $[3n,2k,2d]_2$-code.
Observe $d=n-s.$ The Griesmer bound yields $3s\geq (n-s)(2^{2k-2}-1)/2^{2k-2}.$ Solve for $n/s.$
This shows $n/s\leq 4+3/(4^{k-1}-1).$\\
In order to obtain a lower bound on $n/s$ consider the additive code described by choosing as codelines all the lines
of $PG(2k-1,2).$ We have $n=g(2k), s=g(2k-1)$ and 
$$n/s=g(2k)/g(2k-1)=4\times\frac{(4^k-1)(4^k-2)}{(4^k-2)(4^k-4)}=4+3/(4^{k-1}-1).$$
\end{proof}

\section{A family of optimal additive codes}
\label{shortsection}

\begin{Definition}
\label{minsdef}
For each $k$ let $s_k$ be the minimal $s$ such that $n_k(s)/s=\lambda_k.$
\end{Definition}

Because of the sum construction all multiples of $s_k$ also have this property.
We have $\lambda_{1.5}=7$ and $s_{1.5}=1$ corresponding to a 
$[7,1.5,6]_4$-code (codelines are the lines of the Fano plane). 
Then $\lambda_{2}=5, s_2=1$ corresponding to a $[5,2,4]_4$-code
(codelines are the lines of a spread of $PG(3,2),$ equivalently the projective line $PG(1,4)$).
We have $\lambda_{2.5}=31/7$ with $[31,2.5,24]_4$ as an optimal code, hence $s_{2.5}=7$ and
$\lambda_{3}=21/5=4.2$ because of $PG(2,4),$ hence $s_3=5.$ Finally $PG(3,4)$ shows 
$\lambda_{4}=4+1/21=85/21$ and $s_4=21.$
Certainly $s_{3.5}\leq g(6)=21\times 31.$
As $\lambda_{3.5}=127/31$ it follows that $s_{3.5}$ is a multiple of $31.$
We will see below that we have equality.

\begin{Definition}
\label{coverdef}
An $m$-cover of $PG(l-1,2)$ is a set of lines (the codelines, including multiplicities) such that each point
$P$ is on precisely $m$ codelines.
\end{Definition}

Clearly $PG(l-1,2)$ for even $l$ has a 1-cover (a spread of lines), and a $3$-cover of $PG(l-1,2)$
consists of $2^l-1$ lines.

\begin{Lemma}
\label{coverlemma}
$PG(l-1,2)$ for odd $l$ has a 3-cover.
\end{Lemma}
\begin{proof}
Let $E$ be a Fano subplane. The Blokhuis-Brouwer construction \cite{BB} shows that there is a partial spread
covering the points outside $E.$ Use these as codelines, each with multiplicity $3.$ The 3-cover is  completed
by using the lines of $E$ as codelines.
\end{proof}

\begin{Theorem}
\label{asympGriesmertheorem}
We have $s_k=(4^{k-1}-1)/3$ if $k$ is integer and $s_k=4^{k-1}-1$ if $k$ is not integer.
\end{Theorem}
\begin{proof}
Consider $\lambda_k$ in Theorem \ref{asympGriesmertheorem}. As numerator and denominator have a gcd of $3$
when $k$ is integer whereas they are coprime when $k$ is not integer we have that $(4^{k-1}-1)/3$ divides $s_k$ in the first case,
and $4^{k-1}-1$ divides $s_k$ in the latter case. We claim that a 1-cover and a 3-cover of $PG(2k-1,2)$ respectively
achieve those values. \\
Consider the case of the 1-cover when $k$ is an integer. We have $n=(4^k-1)/3.$ Let $H$ be a hyperplane, $s$
the number of codelines contained in $H,t$ the number of remaining codelines. Then
$$3s+t=\vert H\vert , s+t=n.$$
Solving for $s$ shows $s=s_k=(4^{k-1}-1)/3.$ We confirm that $n/s=\lambda_k.$\\
Let now $k=(2i+1)/2$ and let the codelines form a 3-cover of $PG(2i,2).$ Then $n=2^{2i+1}-1.$
Let $s$ be the number of codelines in hyperplane $H,t$ the number of remaining codelines. Then
$$3s+t=3\vert H\vert , s+t=n.$$
Solving for $s$ shows $s=2^{2i-1}-1=4^{k-1}-1.$ Again $n/s=\lambda_k.$
\end{proof}

The proof of Theorem \ref{asympGriesmertheorem} shows that our codes have constant weight. 
The case of integer $k$ is not of independent interest. We may choose
$PG(k-1,4).$ In the other case $k=(2i+1)/2$ we have codes
$$[s_k\lambda_k,k,s_k(\lambda_k-1)]_4=[4^k-1,k,4^k-4^{k-1}]_4$$
which are optimal in a strong sense. The binary concatenated codes meet the Griesmer bound with equality.
We confirm for example $s_{2.5}=7,$ code $[31,2.5,24]_4$ and
$s_{3.5}=31,$ code $[127,3.5,96]_4.$ In particular $n_{3.5}(31)=127.$

\begin{Theorem}
\label{3covertheorem}
Let $k=l/2$ where $l$ is an odd integer. A multiset ${\cal F}$ of lines in $PG(2k-1,2)$ defines a code
$[4^k-1,k,4^k-4^{k-1}]_4=[2^{2k}-1,k,3\cdot 2^{2k-2}]_4$ if and only if it is a 3-cover.
\end{Theorem}

\begin{proof}
We saw one direction in the proof of Theorem \ref{asympGriesmertheorem}.
Let now ${\cal F}$ describe a code $[2^{2k}-1,k,3\cdot 2^{2k-2}]_4.$ 
Double-counting incident pairs $(g,H)$ where $g$ is a codeline and $H$ a hyperplane, we see that
we have constant-weight codes.
The concatenated code is a binary linear
$[3(2^{2k}-1),2k,6\cdot 2^{2k-2}]_2$-code. Let $G$ be a generator matrix. Clearly the Griesmer bound
is met with equality. Our claim is that each point $P\in PG(l-1,2)$ appears precisely three times as a column of $G.$
Start with the residual code, a $[3(2^{2k-1}-1),2k-1,3\cdot 2^{2k-2}]_2$-code.
Continue this process. It ends with a $[3,1,3]$-code. This shows that point   $(0:\dots :0:1)^T$  appears at least $3$
times as a column of $G.$ As this is true for any point in $PG(l-1,2),$ our claim is proved.
\end{proof}

 \section{A variant}
\label{variantsection}

 \begin{Theorem}
\label{3.5gentheorem}
Let $m\geq 1$ be an integer and let an additive quaternary $(m+0.5)$-dimensional code $C$
be geometrically described by the following set of lines in $PG(2m,2):$
\begin{itemize}
\item A partial spread partitioning the points outside a Fano plane $E,$
\item three different lines of $E.$
\end{itemize}
Then $C$ is a 2-weight $[(2^{2m+1}+1)/3,2m+1,2^{2m-1}]_4$-code with second weight
$w=d+1=2^{2m-1}+1.$
\end{Theorem}
\begin{proof}
Let $n_0$ be the number of lines in the partial spread. Then $n_0=(2^{2m+1}-8)/3$ and our code has
length $n=n_0+3=(2^{2m+1}+1)/3.$
Let $H$ be a hyperplane containing $E$ and let $x$ be the number of lines in the partial spread contained in $H.$
Then $3x+(n_0-x)$ is the number of points of $H$ outside $E,$ which is $2^{2m}-8.$
It follows $x=(2^{2m-1}-8)/3$ and the number of codelines in $H$ is $(2^{2m-1}+1)/3.$
The corresponding codeword weight is $2^{2m-1}=d.$\\
If $H$ does not contain $E$ then $H\cap E$ may or may not be a codeline.
We have $2x+n_0=2^{2m}-4$ and $x=(2^{2m-1}-2)/3.$ The number of codelines in $H$ is
$x$ or $x+1.$ The corresponding codeword weight is $d+1$ or $d.$
\end{proof}
Case $m=2$ yields two $[11,2.5,8]_4$-codes, case $m=3$ yields two $[43,3.5,32]_4$-codes.
In case $m=4$ we obtain $[171,4.5,128]_4$ which is
interesting as it is known that the largest minimum distance of a linear $[171,5]_4$-code is $127.$\\
{\bf Data availability:} All data generated or analyzed during this study are included in this published article.\\
{\bf Declarations}\\
{\bf Conflict of interest} On behalf of all authors, the corresponding author states
that there is no conflict of interest.


\begin{thebibliography}{99}
\bibitem{book2nded} Bierbrauer J.: Introduction to Coding Theory. Second Edition,
Chapman and Hall/CRC Press, (2017)

\bibitem{dim2pt5} Bierbrauer J., Marcugini S., Pambianco F.:
 Optimal additive quaternary codes of low dimension.
  IEEE IT Transactions {\bf 67},  5116-5118 (2021)

  \bibitem{BB} Blokhuis A., Brouwer A.E.:
Small additive quaternary codes. European Journal of Combinatorics {~\bf25}, 161-167 (2004)

\bibitem{quantgeom}  Bierbrauer J., Faina G., Giulietti M., Marcugini S., Pambianco F.:
  The geometry of quantum codes. Innov.  Incidence Geom. {\bf 6}, 53--71 (2009)

\bibitem{CRRS} Calderbank A.R., Rains E.M., Shor P.M., Sloane N.J.A.:
 Quantum error-correction via codes over $GF(4),$.
 IEEE Transactions on Information Theory {\bf44}, 1369-1387 (1998)

\bibitem{DHY} Dong Y., Hu D., Yu S.:
  Breeding quantum error-correcting codes. Phys. Rev. A {\bf 81}, 022322 (2010)

\bibitem{ES2016} Etzion T., Storme L.:
Galois geometries and coding theory. Des. Codes Cryptogr. {\bf 78}, 311–350 (2016)

\bibitem{Giulietti2013} Giulietti M.: The geometry of covering codes: small complete caps and saturating sets in Galois spaces. Surveys in Combinatorics 2013 , (eds. S. R. Blackburn, R. Holloway and M. 13 Wildon) London Math. Soc., Lecture Notes Series {\bf  409}, Cambridge University Press, Cambridge, 51–90  (2013)

 \bibitem{GLLL17} Guo L.B., Liu Y., Lu L.D., Li R.H.: On construction of good quaternary additive codes. ITM Web of Sciences {\bf 12}, 03013 (2017)

 \bibitem{Hirschfeld1985} Hirschfeld J.W.P.: Finite Projective Spaces of Three Dimensions. Oxford Univ. Press, Oxford, (1985)

 \bibitem{Hirschfeld1999} Hirschfeld J.W.P.: Projective Geometries over Finite Fields. 2nd edition,  Oxford Univ. Press, Oxford, (1999)

\bibitem{HirschfeldStorme}Hirschfeld J.W.P., Storme L.: The packing problem in statistics, coding theory and finite projective spaces: Update 2001. Finite Geometries (Proc. 4th Isle of Thorns Conf., July 16-21, 2000),  eds. Blokhuis A., Hirschfeld J.W.P., Jungnickel D., Thas J.A., Develop. Math., 3, Kluwer, Dordrecht, 201–246 (2001)

\bibitem{LandgevStorme} Landjev I., Storme L., Galois geometry and coding theory. Current Research Topics in Galois geometry, eds. De Beule J., Storme L. , Chapter 8, NOVA Academic, New York, {\bf 31}, 187–214 (2011)

 \bibitem{SLKS23} Shi M., Liu L., Kim J.L., Sol\'e P.: Additive complementary dual codes over $GF(4),$.
 Des. Codes Cryptogr. {\bf 91}, 273-284  (2023)


 \end{thebibliography}
\end{document}